\documentclass[a4paper,12pt]{article}
\usepackage{amsfonts,amssymb,latexsym,amsmath}
\usepackage{multicol}
\usepackage[cp1251]{inputenc}
\usepackage[russian]{babel}
\begin{document}
\newcounter{bnomer}
\newcounter{snomer}
\newcounter{diagram}
\setcounter{bnomer}{0} \setcounter{diagram}{0}
\renewcommand{\thesnomer}{\thebnomer.\arabic{snomer}}
\renewcommand{\thebnomer}{\arabic{bnomer}}

\newcommand{\sect}[1]{%
\setcounter{snomer}{0} \refstepcounter{bnomer}
\begin{center}\large{\textbf{\S \thebnomer.{ #1}}}\end{center}}

\newcommand{\thenv}[2]{%
\refstepcounter{snomer}
\par\addvspace{\medskipamount}\textbf{{#1} \thesnomer.}
{#2}\par\addvspace{\medskipamount}}


\renewcommand{\refname}{References}

\date{}
\title{The algebra of invariants for the adjoint action of the
unitriangular group}
\author{Victoria Sevostyanova\thanks{Partially
supported by Israel Scientific Foundation grant 797/14.}}

\maketitle

\begin{center}
\parbox[b]{330pt}{\small\textsc{Abstract.} In this paper the
algebra of invariants for the adjoint action of the unitriangular
group in the nilradical of a parabolic subalgebra is studied. We
prove that the algebra of invariants is finitely generated.}
\end{center}

\sect{Introduction}

Let $G$ be the general linear group $\mathrm{GL}(n,K)$ over an
algebraically closed field $K$ of characteristic zero. Let $B$ ($N$,
respectively) be its Borel (maximal unipotent, respectively)
subgroup, which consists of upper triangular matrices with nonzero
(unit, respectively) elements on the diagonal. We fix a parabolic
subgroup $P\supset B$. Let $\mathfrak{p}$, $\mathfrak{b}$ and
$\mathfrak{n}$ be the Lie subalgebras in $\mathfrak{gl}(n,K)$
correspon\-ding to $P$, $B$ and $N$, respectively. We represent
$\mathfrak{p}=\mathfrak{r}\oplus\mathfrak{m}$ as the direct sum of
the nilradical $\mathfrak{m}$ and a block diagonal subalgebra
$\mathfrak{r}$ with sizes of blocks $(n_1,\ldots, n_s)$. The
subalgebra $\mathfrak{m}$ is invariant relative to the adjoint
action of the group $P$:
$$\mbox{for any }g\in P\mbox{ we have
}x\in\mathfrak{m}\mapsto\mathrm{Ad}_gx=gxg^{-1}.$$
Therefore $\mathfrak{m}$ is invariant relative to the adjoint action
of the subgroups $B$ and $N$. We extend this action to the
representation in the algebra $K[\mathfrak{m}]$ and in the field
$K(\mathfrak{m})$:
$$\mbox{for any }g\in P\mbox{ we have
}f(x)\in K[\mathfrak{m}]\mapsto f(\mathrm{Ad}_{g^{-1}}x).$$

The complete description of the field of invariants
$K(\mathfrak{m})^N$ for any parabolic subalgebra is a result
of~\cite{S1}. In this paper a notion of an extended base is
introduced. Elements of the extended base correspond to a set of
algebraically independent $N$-invariants. These invariants generate
the field of invariants $K(\mathfrak{m})^N$. Further in the
paper~\cite{S2} the structure of the algebra of invariants
$K[\mathfrak{m}]^N$ is considered. If the sizes of diagonal blocks
are $(2,k,2)$, $k>2$, or $(1,2,2,1)$, then the invariants
constructed on the extended base do not generate the algebra of
invariants and the algebra of invariants is not free. Besides, the
additional invariants in both cases are constructed, which together
with the main list of the invariants constructed on the extended
base generate the algebra of invariants. Also, the relations between
these invariants are provided.

The aim of this paper is to prove that the algebra of invariants
$K[\mathfrak{m}]^N$ is finitely generated. We show this as follows.
Let $P=L\ltimes U$, where $L$ is the Levi subgroup and $U$ is the
unipotent radical. Then $N=U_L\ltimes U$, where $U_L$ is the maximal
unipotent subgroup of $L$. One has
$$K[\mathfrak{m}]^N=K\Big[K[\mathfrak{m}]^U\Big]^{U_L}.$$
In this paper we show that the algebra of invariants
$K[\mathfrak{m}]^U$ is a finitely generated, free algebra and we
present its generating invariants. Then by Khadzhiev's theorem (see
Theorem~\ref{Th-Khadzhiev}), we get our main result:

\medskip
\thenv{Theorem}{\emph{The algebra of invariants $K[\mathfrak{m}]^N$
is finitely generated}.}


\sect{Main statements and definitions}

We begin with definitions. Let
$\mathfrak{b}=\mathfrak{n}\oplus\mathfrak{h}$ be a triangular
decomposition. Let $\Delta$ be the root system relative to
$\mathfrak{h}$ and let $\Delta^{\!+}$ be the set of positive roots.
Let $\{\varepsilon_i\}_{i=1}^{n}$ be the standard basis of
$\mathbb{C}^n$. Every positive root $\gamma$ in $\mathfrak{gl}(n,K)$
can be represented as $\gamma=\varepsilon_i-\varepsilon_j$,
$1\leqslant i<j\leqslant n$ (see \cite{GG}). We identify a root
$\gamma$ with the pair $(i,j)$ and the set of the positive roots
$\Delta^{\!+}$ with the set of pairs $(i,j)$, $i<j$. The system of
positive roots $\Delta^{\!+}_\mathfrak{r}$ of the reductive
subalgebra $\mathfrak{r}$ is a subsystem in $\Delta^{\!+}$.

Let $\{E_{i,j}:~i<j\}$ be the standard basis in $\mathfrak{n}$. Let
$E_\gamma$ denote the basis element $E_{i,j}$, where $\gamma=(i,j)$.

Let $M$ be a subset of $\Delta^{\!+}$ corresponding to
$\mathfrak{m}$ that is $$\mathfrak{m}=\bigoplus_{\gamma\in
M}E_{\gamma}.$$ We identify the algebra $K[\mathfrak{m}]$ with the
polynomial algebra in variables $x_{i,j}$, $(i,j)\in M$.

We define a relation in $\Delta^{\!+}$ such that $\gamma'>\gamma$
whenever $\gamma'-\gamma\in\Delta^{\!+}$. Note that the relation $>$
is not an order relation.

The roots $\gamma$ and $\gamma'$ are called \emph{comparable}, if
either $\gamma'>\gamma$ or $\gamma>\gamma'$.

We will introduce a subset $S$ in the set of positive roots such
that every root from this subset corresponds to some $N$-invariant.

\thenv{Definition}{A subset $S$ in $M$ is called a \emph{base} if
the elements in $S$ are not pairwise comparable and for any
$\gamma\in M\setminus S$ there exists $\xi\in S$ such that
$\gamma>\xi$.}

Let us show that the base exists. We need the following

\thenv{Definition\label{Def-minimal}}{Let $A$ be a subset in $M$. We
say that $\gamma$ is a \emph{minimal element} in $A$ if there is no
$\xi\in A$ such that $\xi<\gamma$.}

For a given parabolic subgroup we will construct a diagram in the
form of a square array. The cell of the diagram corresponding to a
root of $S$ is labeled by the symbol $\otimes$. Symbols $\times$
will be explained below.

\thenv{Example\label{Ex-2132}}{Diagram 1 represents the parabolic
subalgebra with sizes of its diagonal blocks $(2,1,3,2)$. In this
case minimal elements in $M$ are $(2,3)$, $(3,4)$ and $(6,7)$.}
\begin{center}\refstepcounter{diagram}
{\begin{tabular}{|p{0.1cm}|p{0.1cm}|p{0.1cm}|p{0.1cm}|p{0.1cm}|p{0.1cm}|p{0.1cm}|p{0.1cm}|c}
\multicolumn{2}{l}{{\small 1\quad2\!\!}}&\multicolumn{2}{l}{{\small
3\quad4\!\!}}&\multicolumn{2}{l}{{\small 5\quad 6\!\!}}&
\multicolumn{2}{l}{{\small 7\quad 8\!\!}}\\
\cline{1-8} \multicolumn{2}{|l|}{1}&&&$\!\otimes$&&&&{\small 1}\\
\cline{3-8} \multicolumn{2}{|r|}{1}&$\!\otimes$&&&&&&{\small 2}\\
\cline{1-8} \multicolumn{2}{|c|}{}&1&$\!\otimes$&&&&&{\small 3}\\
\cline{3-8} \multicolumn{3}{|c|}{}&\multicolumn{3}{|l|}{1}&$\!\times$&$\!\times$&{\small 4}\\
\cline{7-8} \multicolumn{3}{|c|}{}&\multicolumn{3}{|c|}{1}&$\!\times$&$\!\otimes$&{\small 5}\\
\cline{7-8} \multicolumn{3}{|c|}{}&\multicolumn{3}{|r|}{1}&$\!\otimes$&&{\small 6}\\
\cline{4-8} \multicolumn{6}{|c|}{}&\multicolumn{2}{|l|}{1}&{\small 7}\\
\multicolumn{6}{|c|}{}&\multicolumn{2}{|r|}{1}&{\small 8}\\
\cline{1-8} \multicolumn{8}{c}{Diagram \arabic{diagram}}\\
\end{tabular}}
\end{center}

We construct the base $S$ by the following algorithm.

\textsc{Step 1.} Put $M_0=M$ and $i=1$. Let $S_1$ be the set of
minimal elements in $M_0$.

\textsc{Step 2.} Put $M_i=M_{i-1}\setminus\big\{S_i\cup\{\gamma\in
M_{i-1}:\exists\xi\in S_i, \xi<\gamma\}\big\}$. Let $S_i$ be the set
of minimal elements of $M_{i-1}$. Increase $i$ by 1 and repeat Step
2 until $M_i$ is empty.

Denote $S=S_1\cup S_2\cup\ldots$ The base $S$ is unique.

We have $S_1=\{(2,3),(3,4),(6,7)\}$ and $S_2=\{(1,5),(5,8)\}$ in
Example~\ref{Ex-2132}.

\medskip
Let $(r_1,r_2,\ldots,r_s)$ be the sizes of the diagonal blocks in
$\mathfrak{r}$. Put $$R_k=\displaystyle\sum_{i=1}^{k}r_i.$$

Let us present $N$-invariant corresponding to a root of the base.
Consider the formal matrix $\mathbb{X}$ of variables
$$\left(\mathbb{X}\right)_{i,j}=\left\{
\begin{array}{ll}
x_{i,j}&\mbox{if }(i,j)\in M;\\
0&\mbox{otherwise.}
\end{array}\right.$$
The matrix $\mathbb{X}$ can be represented as a block matrix
$$\mathbb{X}=\left(\begin{array}{ccccc}
0&X_{1,2}&X_{1,3}&\ldots&X_{1,s}\\
0&0&X_{2,3}&\ldots&X_{1,s}\\
\ldots&\ldots&\ldots&\ldots&\ldots\\
0&0&0&\ldots&X_{s-1,s}\\
0&0&0&\ldots&0\\
\end{array}\right),$$
where the size of $X_{i,j}$ is $r_i\times r_j$,
\begin{equation}
X_{i,j}=\left(\begin{array}{cccc}
x_{R_{i-1}+1,R_{j-1}+1}&x_{R_{i-1}+1,R_{j-1}+2}&\ldots&x_{R_{i-1}+1,R_{j}}\\
x_{R_{i-1}+1,R_{j-1}+2}&x_{R_{i-1}+2,R_{j-1}+2}&\ldots&x_{R_{i-1}+2,R_{j}}\\
\ldots&\ldots&\ldots&\ldots\\
x_{R_{i},R_{j-1}+1}&x_{R_{i},R_{j-1}+2}&\ldots&x_{R_{i},R_{j}}\\
\end{array}\right).\label{X_ij}
\end{equation}

\thenv{Lemma\label{Lemma1}}{\emph{The roots corresponding to the
antidiagonal elements in} $X_{i,i+1}$ (\emph{from the lower left
element towards right upper direction}) \emph{are in the base}.}

Thus the roots of the base in the blocks $X_{i,i+1}$ are as follows.
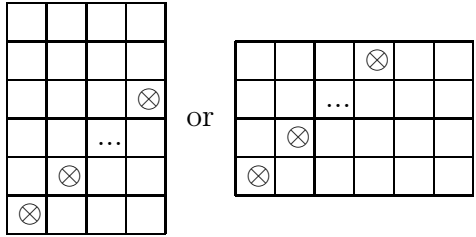
\begin{center}
{\begin{tabular}{|p{0.1cm}|p{0.1cm}|p{0.1cm}|p{0.1cm}|}
\cline{1-4} &&&\\
\cline{1-4} &&&\\
\cline{1-4} &&&$\!\otimes$\\
\cline{1-4} &&$\!...$&\\
\cline{1-4} &$\!\otimes$&&\\
\cline{1-4} $\!\otimes$&&&\\
\cline{1-4}
\end{tabular}\ \ \mbox{or }
\begin{tabular}{|p{0.1cm}|p{0.1cm}|p{0.1cm}|p{0.1cm}|p{0.1cm}|p{0.1cm}|}
\cline{1-6} &&&$\!\otimes$&&\\
\cline{1-6} &&$\!...$&&&\\
\cline{1-6} &$\!\otimes$&&&&\\
\cline{1-6} $\!\otimes$&&&&&\\
\cline{1-6}
\end{tabular}}
\end{center}

\textsc{Proof.} By definition~\ref{Def-minimal} for any $i$ the root
$(R_i,R_i+1)$ is minimal. Therefore $M\setminus M_1$ contains roots
corresponding to all cells in the row $R_i$ and the column $R_i+1$.
Hence $(R_i-1,R_i+2)\in S_2$ if $r_i,r_{i+1}>1$ and all roots of $M$
in the rows $R_i$, $R_i-1$ and in the columns $R_i+1$, $R_i+2$
belong to $M\setminus M_2$. Hence $(R_i-2,R_i+3)\in S_3$ if
$r_i,r_{i+1}>2$ etc.~$\Box$

There are roots in $S$ such that these roots do not correspond to
elements of the secondary diagonal in $X_{i,i+1}$, for example
$(1,5)$ in Example~\ref{Ex-2132}.

\medskip
For any root $\gamma=(a,b)\in M$ let $S_\gamma=\{(i,j)\in
S:i>a,j<b\}$. Let $S_\gamma=\{(i_1,j_1),\ldots,(i_k,j_k)\}$. Note
that if $\gamma$ is minimal in $M$, then $S_{\gamma}=\emptyset$.
Denote by $M_\gamma$ a minor $\mathbb{X}_I^J$ of the matrix
$\mathbb{X}$ with ordered systems of rows $I$ and columns $J$, where
$$I=\mathrm{ord}\{a,i_1,\ldots,i_k\},\quad J=\mathrm{ord}\{j_1,\ldots,j_k, b\}.$$

\thenv{Example}{Let us continue Example~\ref{Ex-2132}. For the root
$(1,6)$ we have $S_{(1,6)}=\{(2,3),(3,4)\}$, $I=\{1,2,3\}$,
$J=\{3,4,6\}$, and
$$M_{(1,6)}=\left|\begin{array}{ccc}
x_{1,3}&x_{1,4}&x_{1,6}\\
x_{2,3}&x_{2,4}&x_{2,6}\\
0&x_{3,4}&x_{3,6}\\
\end{array}\right|.$$
All minors $M_{\xi}$ for $\xi\in S$ are following
$$M_{(2,3)}=x_{2,3},\ M_{(3,4)}=x_{3,4},\ M_{(6,7)}=x_{6,7},$$$$
M_{(5,8)}=\left|\begin{array}{cc}
x_{5,7}&x_{5,8}\\
x_{6,7}&x_{6,8}\\
\end{array}\right|,\
M_{(1,5)}=\left|\begin{array}{ccc}
x_{1,3}&x_{1,4}&x_{1,5}\\
x_{2,3}&x_{2,4}&x_{2,5}\\
0&x_{3,4}&x_{3,5}\\
\end{array}\right|.$$}

\thenv{Lemma\label{L-M-inv}}{\emph{For any $\xi\in S$ the minor
$M_{\xi}$ is $N$-invariant}.}

\thenv{Notation\label{Notation}}{The group $N$ is generated by the
one-parameter subgroups
$$g_{i,j}(t)=I+tE_{i,j},\mbox{ where }1\leqslant i<j\leqslant n$$
and $I$ is the identity matrix. The adjoint action of any
$g_{i,j}(t)$ makes the following transformations of a matrix:
\begin{itemize}
\item[1)] the $j$th row multiplied by $t$ is added to the $i$th row,
\item[2)] the $i$th column multiplied by $-t$ is added to the $j$th
column, i.e. for a variable $x_{a,b}$ we have
$$\mathrm{Ad}_{g_{i,j}^{-1}(t)}x_{a,b}=\left\{\begin{array}{ll}
x_{a,b}+tx_{j,b}&\mbox{if }a=i;\\
x_{a,b}-tx_{a,i}&\mbox{if }b=j;\\
x_{a,b}&\mbox{otherwise}.
\end{array}\right.$$
\end{itemize}}

\textsc{Proof.} By the notation it is sufficient to prove that for
any $\xi=(k,m)\in S$ the minor $M_{\xi}$ is invariant under the
adjoint action of $g_{i,j}(t)$ for any $i<j$. If $i<k$, then the
$i$th row does not belong to the minor $M_{\xi}$ and the adding of
the $j$th row to the $i$th row leaves $M_{\xi}$ unchanged. Let
$M_{\xi}=\mathbb{X}_I^J$ for some collections of rows $I$ and
columns $J$. If $i\geqslant k$, then since the numbers in $I$ are
consecutive, the number of any nonzero row $j$ at the intersection
with columns $J$ belongs to $I$. Then the adding of the $j$th row to
the $i$th row leaves $M_{\xi}$ unchanged again. Using the similar
reasoning for columns, we get that $M_{\xi}$ is
$N$-invariant.~$\Box$

\medskip
The set $\{M_{\xi},\ \xi\in S\}$ does not generate all the
$N$-invariants. There is the other series of $N$-invariants. To
present it we need

\thenv{Definition}{An ordered set of positive roots
$$\{\varepsilon_{i_1}-\varepsilon_{j_1},
\varepsilon_{i_2}-\varepsilon_{j_2},\ldots,
\varepsilon_{i_s}-\varepsilon_{j_s}\}$$ is called a \emph{chain} if
$j_1=i_2,j_2=i_3,\ldots,j_{s-1}=i_s$.}

\thenv{Definition}{We say that two roots $\xi,\xi'\in S$ form an
\emph{admissible pair} $q=(\xi,\xi')$ if there exists $\alpha_q$ in
the set $\Delta^{\!+}_\mathfrak{r}$ corresponding to the reductive
part $\mathfrak{r}$ such that the ordered set of roots
$\{\xi,\alpha_q,\xi'\}$ is a chain. In other words, roots
$\xi=\varepsilon_i-\varepsilon_j$ and
$\xi'=\varepsilon_k-\varepsilon_l$ are an admissible pair if
$\alpha_q=\varepsilon_j-\varepsilon_k\in\Delta^{\!+}_\mathfrak{r}$.
Note that the root $\alpha_q$ is uniquely determined by $q$.}

\thenv{Example\label{Ex-2132-admissible_pair}}{In the case of
Diagram 1 we have three admissible pairs $q_1=(\xi_1,\xi_3)$,
$q_2=(\xi_2,\xi_3)$, $q_3=(\xi_1,\xi_4)$, where $\xi_1=(2,3)$,
$\xi_2=(1,5)$, $\xi_3=(6,7)$, and $\xi_4=(5,8)$.}

Let the set $Q:=Q(\mathfrak{p})$ consist of admissible pairs. For
every admissible pair $q=(\xi,\xi')$ we construct a positive root
$\varphi_q=\alpha_q+\xi'$, where $\{\xi,\alpha_q,\xi'\}$ is a chain.
Consider the subset $\Phi=\{\varphi_q:~q\in Q\}$ in the set of
positive roots. The cell of the diagram corresponding to a root of
$\Phi$ is labeled by $\times$.

\thenv{Example}{The roots of $\Phi$ for the admissible pairs in
Example~\ref{Ex-2132-admissible_pair} are $\varphi_{q_1}=(4,7)$,
$\varphi_{q_2}=(5,7)$, $\varphi_{q_3}=(4,8)$.}

Now we are ready to present the $N$-invariant corresponding to a
root $\varphi\in\Phi$.

Let admissible pair $q=(\xi,\xi')$ correspond to $\varphi_q\in\Phi$.
We construct the polynomial
\begin{equation}
L_{\varphi_q}=\sum_{\scriptstyle\alpha_1,\alpha_2\in\Delta^{\!+}_\mathfrak{r}\cup\{0\}
\atop\scriptstyle\alpha_1+\alpha_2=\alpha_q}
M_{\xi+\alpha_1}M_{\alpha_2+\xi'}.\label{L_q}
\end{equation}

\thenv{Example}{Continuing the previous example, we have
$$L_{(4,7)}=x_{3,4}x_{4,7}+x_{3,5}x_{5,7}+x_{3,6}x_{6,7},$$
$$L_{(4,8)}=x_{3,4}\left|\begin{array}{cc}
x_{4,7}&x_{4,8}\\
x_{6,7}&x_{6,8}\\
\end{array}\right|+x_{3,5}\left|\begin{array}{cc}
x_{5,7}&x_{5,8}\\
x_{6,7}&x_{6,8}\\
\end{array}\right|,$$
$$L_{(5,7)}=\left|\begin{array}{ccc}
x_{1,3}&x_{1,4}&x_{1,5}\\
x_{2,3}&x_{2,4}&x_{2,5}\\
0&x_{3,4}&x_{3,5}\\
\end{array}\right|x_{5,7}+\left|\begin{array}{ccc}
x_{1,3}&x_{1,4}&x_{1,6}\\
x_{2,3}&x_{2,4}&x_{2,6}\\
0&x_{3,4}&x_{3,6}\\
\end{array}\right|x_{6,7}.$$}

\thenv{Lemma}{\emph{The polynomial $L_{\varphi}$ is $N$-invariant
for any} $\varphi=\varphi_q\in\Phi$, $q=(\xi,\xi')$.}

\textsc{Proof.} By Notation~\ref{Notation} it is sufficient to check
the action of $g_{i,j}(t)$. Let $\xi=(a,b)$,~$\xi'=(a',b')$. Using
the definition of admissible pair, we have $a<b<a'<b'$,
$\alpha_q=(b,a')\in\Delta_{\mathfrak{r}}^{\!+}$, and
$\varphi=(b,b')$. If $i<b$ or $j>a'$, then using the same arguments
as in the proof of the invariance of $M_{\xi}$ for $\xi\in S$, we
have that the minors of the right part of (\ref{L_q}) are
$g_{i,j}(t)$-invariant.

Let $b\leqslant i<j\leqslant a'$. Denote $\gamma_1=(b,i)$,
$\gamma_2=(j,a')$, $\beta=(i,j)$; then
$\alpha_q=\gamma_1+\beta+\gamma_2$ and
$\gamma_1+\beta,\beta+\gamma_2\in\Delta_{\mathfrak{r}}^{\!+}\cup\{0\}$.
We have
\begin{equation}
\left\{\begin{array}{l} T_{g_{i,j}(t)}M_{\xi+\gamma_1+\beta}=
M_{\xi+\gamma_1+\beta}+tM_{\xi+\gamma_1},\\
T_{g_{i,j}(t)}M_{\beta+\gamma_2+\xi'}=
T_{\beta+\gamma_2+\xi'}-tM_{\gamma_2+\xi'}.\\
\end{array}\right.\label{T_g(M_xi)}
\end{equation}
The other minors of (\ref{L_q}) are invariant under the action of
$g_{i,j}(t)$. Combining (\ref{L_q}) and (\ref{T_g(M_xi)}), we get
$$\left(T_{g_{i,j}(t)}L_{\varphi}\right)-L_{\varphi}=M_{\xi+\gamma_1}
\left(M_{\beta+\gamma_2+\xi'}-t M_{\gamma_2+\xi'}\right)+$$
$$\left(M_{\xi+\gamma_1+\beta}+t
M_{\xi+\gamma_1}\right)M_{\gamma_2+\xi'}-
M_{\xi+\gamma_1}M_{\beta+\gamma_2+\xi'}-
M_{\xi+\gamma_1+\beta}M_{\gamma_2+\xi'}=0.~\Box$$

Thus we proved the first part of

\thenv{Theorem\label{M-L_independ}}{\emph{For an arbitrary parabolic
subalgebra}, \emph{the system of polynomials}
\begin{equation}
\{M_\xi,~\xi\in S,~L_{\varphi},~\varphi\in\Phi,\}\label{system-M-L}
\end{equation}
\emph{is contained in $K[\mathfrak{m}]^N$ and is algebraically
independent over $K$}.}

To show the algebraic independence, consider the restriction
homo\-mor\-phism $f\mapsto f|_\mathcal{Y}$, where
$$\mathcal{Y}=\left\{\sum_{\xi\in S\cup\Phi}c_{\xi}E_\xi:
c_{\xi}\neq0\ \forall\xi\in S\cup\Phi\right\},$$ from
$K[\mathfrak{m}]$ to the polynomial algebra $K[\mathcal{Y}]$ of
$x_\xi$, $\xi\in S$, and of $x_\varphi$, $\varphi\in\Phi$. Direct
calculations show that the system of the images
$$\left\{M_\xi|_{\mathcal{Y}},\xi\in S,L_\varphi|_{\mathcal{Y}},
\varphi\in\Phi\right\}$$ is algebraically independent over $K$.
Therefore, the system~(\ref{system-M-L}) is algebraically
independent over $K$ (see details in~\cite{PS}).

\thenv{Definition}{The set $S\cup\Phi$ is called an \emph{extended
base}.}

\thenv{Definition}{The matrices of $\mathcal{Y}$ are called
\emph{canonical}.}

By~\cite{S1} one has the following theorems.

\thenv{Theorem\label{Exist_of_representative}}{\emph{There exists a
nonempty Zariski-open subset $W\subset\mathfrak{m}$ such that the
$N$-orbit of any $x\in W$ intersects $\mathcal{Y}$ at a unique
point}.}

\thenv{Theorem\label{Th_invariant_field}}{\emph{The field of
invariants $K(\mathfrak{m})^N$ is the field of rational functions
of} $M_\xi$, $\xi\in S$, \emph{and} $L_{\varphi}$,
$\varphi\in\Phi$.}


\sect{Invariants of the unipotent subgroup \\in the Levi
decomposition of $P$}

Let us consider the decomposition of a parabolic group $P$ into the
semi\-direct product of the Levi subgroup $L$ and the unipotent
radical $U$. Let $U_L$ be the maximal unipotent subgroup in the Levi
group $L$. One has $N=U_L\ltimes U$. The aim is to describe the
algebra of invariants $K[\mathfrak{m}]^{U}$.

As above, we will introduce some subset $T\subset\Delta^{\!+}$ and
construct a correspon\-ding invariant $N_{\xi}\in
K[\mathfrak{m}]^{U}$ for every root $\xi\in T$.

\thenv{Definition}{A root $\xi\in\Delta^{\!+}$ belongs to \emph{a
broad base} $T\subset\Delta^{\!+}$ if one of the following
conditions holds:
\begin{itemize}
\item[1)] the root $\xi$ belongs to $S$;
\item[2)] there exists a root $\gamma\in S$ such that
$\xi>\gamma$ and the variables $x_\xi$ and $x_\gamma$ are located in
the same block $X_{i,j}$.
\end{itemize}}

\thenv{Example\label{Ex-2132-T}}{The diagram presents roots of the
broad base $T$ for the diagonal blocks $(2,1,3,2)$. The cells of the
diagram corresponding to roots of $S$ (resp. $T\setminus S$) are
labeled by the symbol $\otimes$ (resp.~$\boxtimes$).}
\begin{center}\refstepcounter{diagram}
{\begin{tabular}{|p{0.1cm}|p{0.1cm}|p{0.1cm}|p{0.1cm}|p{0.1cm}|p{0.1cm}|p{0.1cm}|p{0.1cm}|c}
\multicolumn{2}{l}{{\small 1\quad 2\!\!}}&\multicolumn{2}{l}{{\small
3\quad 4\!\!}}&\multicolumn{2}{l}{{\small 5\quad 6\!\!}}&
\multicolumn{2}{l}{{\small 7\quad 8\!\!}}\\
\cline{1-8} \multicolumn{2}{|l|}{1}&$\!\boxtimes$&&$\!\otimes$&$\!\boxtimes$&&&{\small 1}\\
\cline{3-8} \multicolumn{2}{|r|}{1}&$\!\otimes$&&&&&&{\small 2}\\
\cline{1-8} \multicolumn{2}{|c|}{}&1&$\!\otimes$&$\!\boxtimes$&$\!\boxtimes$&&&{\small 3}\\
\cline{3-8} \multicolumn{3}{|c|}{}&\multicolumn{3}{|l|}{1}&$\!\boxtimes$&$\!\boxtimes$&{\small 4}\\
\cline{7-8} \multicolumn{3}{|c|}{}&\multicolumn{3}{|c|}{1}&$\!\boxtimes$&$\!\otimes$&{\small 5}\\
\cline{7-8} \multicolumn{3}{|c|}{}&\multicolumn{3}{|r|}{1}&$\!\otimes$&$\!\boxtimes$&{\small 6}\\
\cline{4-8} \multicolumn{6}{|c|}{}&\multicolumn{2}{|l|}{1}&{\small 7}\\
\multicolumn{6}{|c|}{}&\multicolumn{2}{|r|}{1}&{\small 8}\\
\cline{1-8} \multicolumn{8}{c}{Diagram \arabic{diagram}}\\
\end{tabular}}
\end{center}

Let $M'=\big\{\xi\in M:E_{\xi}\in\mathfrak{m}^2\big\}.$ In other
words, if an element corresponding to a root $\xi\in M$ does not
belong to blocks $X_{k,k+1}$ for any $k$, then $\xi\in M'$. We have
$$\sum_{\xi\in M'}x_{\xi}E_{\xi}=\left(\begin{array}{ccccc}
0&0&X_{1,3}&\ldots&X_{1,s}\\
\ldots&\ldots&\ldots&\ldots&\ldots\\
0&0&0&\ldots&X_{s-2,s}\\
0&0&0&\ldots&0\\
0&0&0&\ldots&0\\
\end{array}\right),$$
$$\sum_{\xi\in M\setminus M'}x_{\xi}E_{\xi}=\left(\begin{array}{ccccc}
0&X_{1,2}&0&\ldots&0\\
0&0&X_{1,3}&\ldots&0\\
\ldots&\ldots&\ldots&\ldots&\ldots\\
0&0&0&\ldots&X_{s-1,s}\\
0&0&0&\ldots&0\\
\end{array}\right),$$
where $X_{i,j}$ is the block~(\ref{X_ij}).

If $\xi\in M\setminus M'$, then $x_{\xi}$ is in some block
$X_{k,k+1}$. We have $\xi\in S$ or using Lemma~\ref{Lemma1}, there
is $\gamma\in S$ such that $\xi$ is to the right or above $\gamma$.
In both cases $\xi\in T$. Therefore $M\setminus M'\subset T$.
Example~\ref{Ex-2132-T} shows that $M\setminus M'\neq T$ in general
case.

For $\xi=(i,j)\in T$ let $N_{\xi}\in K[\mathfrak{m}]$ be defined as
follows
$$N_{\xi}=\left\{\begin{array}{ll}
x_{i,j}&\mbox{if }\xi\in M\setminus M';\\
M_{\xi}&\mbox{if }\xi\in M'.
\end{array}\right.$$

\thenv{Example}{Let us write all $U$-invariants $N_{\xi}$ in the
case $(2,1,3,2)$ for $\xi\in T\cap M'$.
$$N_{(1,5)}=\left|\begin{array}{ccc}
x_{1,3}&x_{1,4}&x_{1,5}\\
x_{2,3}&x_{2,4}&x_{2,5}\\
0&x_{3,4}&x_{3,5}\\
\end{array}\right|,\quad
N_{(1,6)}=\left|\begin{array}{ccc}
x_{1,3}&x_{1,4}&x_{1,6}\\
x_{2,3}&x_{2,4}&x_{2,6}\\
0&x_{3,4}&x_{3,6}\\
\end{array}\right|.$$}

\thenv{Lemma\label{Th:N_xi-invariant}}{\emph{The minor $N_{\xi}$ is
invariant under the adjoint action of the unipotent group} $U$
\emph{for any} $\xi\in T$.}

\textsc{Proof.} The group $U$ is generated by the one-parameter
subgroups (see Notation~\ref{Notation})
$$g_{i,j}(t)=I+tE_{i,j},\mbox{ where }(i,j)\in M.$$

There are two cases of a root $\xi\in T$. The first case is $\xi\in
M\setminus M'$ and the second one is $\xi\in M'\cap T$.
\begin{enumerate}
\item
Suppose $\xi=(a,b)\in T$ belongs to the set $M\setminus M'$; then
$N_{\xi}=x_{a,b}$ and there is some $k$ such that the variable
$x_{a,b}$ is in the block $X_{k,k+1}$. Using the
Notion~\ref{Notation}, for $t\neq0$ we have
$\mathrm{Ad}_{g_{i,j}^{-1}(t)}x_{a,b}\neq x_{a,b}$ if $a=i$ and
$x_{j,b}$ is in $X_{k,k+1}$ or $b=j$ and $x_{a,i}$ is in
$X_{k,k+1}$. In both cases the root $(i,j)$ belongs to
$\Delta^{\!+}_{\mathfrak{r}}$. Therefore $(i,j)\not\in M$ and
$g_{i,j}(t)\not\in U$.
Hence $x_{a,b}$ is an $U$-invariant.

\item
If the root $\xi=(a,b)\in T$ does not belong to $M\setminus M'$,
then by definition of $T$, there exists a root $\gamma\in S$ such
that $\gamma=(i,b)$, $i>a$, or $\gamma=(a,j)$, $j<b$, and $x_\xi$
and $x_\gamma$ are in the same block $X_{l,m}$, $l<m+1$. Suppose
$\gamma=(i,b)$, $i>a$. The case $\gamma=(a,j)$ is similar. Let
$M_{\gamma}=\mathbb{X}_I^J$ be a minor of order $k$ of the formal
matrix with rows $I=\{i,i+1,\ldots,i+k-1\}$ and columns
$J=\{b-k+1,b-k+2,\ldots,b\}$, then $N_{\xi}=\mathbb{X}_{I'}^J$,
where $I'=\{a,i+1,\ldots,i+k-1\}$. Note that all rows of $N_{\xi}$
except the row $a$ and all columns are consecutive. Since a minor is
not changed by addition to a row (resp. column) any other its row
(resp. column), the adjoint action of $g_{u,v}(t)$ can change
$\mathbb{X}_{I'}^J$ if $u=a$ and $v\leqslant i$. Let
$\mathrm{Ad}_{g_{u,v}^{-1}(t)}N_{\xi}\neq N_{\xi}$ for $t\neq0$.
Since $x_\xi$ and $x_\gamma$ are in the same block $X_{l,m}$ and
$u=a$ and $v\leqslant i$, then $x_{(u,b)}$ and $x_{(v,b)}$ are in
the same block $X_{l,m}$. Hence
$(u,v)\in\Delta_{\mathfrak{r}}^{\!+}$ and $g_{u,v}(t)\not\in U$. So
we have that $N_{\xi}$ is an $U$-invariant.~$\Box$
\end{enumerate}

\thenv{Definition}{\emph{The remoteness} of a root $\gamma\in M$ is
called the maximum number $s$ of roots $\gamma_i$ in $M$ such that
$\gamma=\gamma_1>\gamma_2>\ldots>\gamma_s$.}

\thenv{Example}{The remoteness of the root $(1,6)$ in
Example~\ref{Ex-2132} equals~5, we have
$$(1,6)>(1,5)>(2,5)>(3,5)>(3,4).$$}

\thenv{Lemma\label{N_independ}}{\emph{The system of polynomials
$\{N_{\xi}, \ \xi\in T\}$ is algebraically independent over $K$}.}

\textsc{Proof.} Assume the converse, namely that the system
$\{N_{\xi}, \ \xi\in T\}$ is alge\-braically dependent. Hence there
is a polynomial $f$ such that for some $\xi_1,\ldots,\xi_k$ we have
$$f(N_{\xi_1},N_{\xi_2},\ldots,N_{\xi_k})=0.$$
Suppose that the degree of the polynomial $f$ is minimal. Let
$\xi_1$ be a root with the maximal remoteness. If $\xi\in T$ has a
$k$th remoteness, then all roots $\gamma\neq\xi$ for variables
$x_{\gamma}$ in the polynomial $N_{\xi}$ have a remoteness smaller
than $\xi$. The variable $x_{\xi}$ is in the first row and the last
column of the minor $N_{\xi}$. Let us expand $N_{\xi}$ according to
the first row. We have $N_{\xi}=ax_{\xi}+b$ for some polynomials $a$
and $b$ and all variables in $a$ and $b$ correspond to the roots
with less remoteness than the remoteness of $\xi$. Then the variable
$x_{\xi_1}$ is included into the single minor $N_{\xi_1}$.

We have
$$0=f(N_{\xi_1},\ldots,N_{\xi_k})=$$
$$=f_m(N_{\xi_2},\ldots,N_{\xi_k})N_{\xi_1}^m+
f_{m-1}(N_{\xi_2},\ldots,N_{\xi_k})N_{\xi_1}^{m-1}+\ldots+
f_0(N_{\xi_2},\ldots,N_{\xi_k}).$$ Since $N_{\xi_1}=ax_{\xi_1}+b$
and $a\not\equiv0$, we conclude that the coefficient of the highest
power for the variable $x_{\xi_1}$ is
$f_m(N_{\xi_2},\ldots,N_{\xi_k})a^m$. Therefore
$$f_m(N_{\xi_2},\ldots,N_{\xi_k})=0.$$
This contradicts the minimality of $f$ and completes the
proof.~$\Box$


\sect{The algebra of $U$-invariants}

Let $\mathcal{Z}=\left\{\displaystyle\sum_{\xi\in
T}c_{\xi}E_{\xi}:c_{\xi}\in K\right\}$.

\thenv{Proposition\label{Exist_of_representative_for_U}}{\emph{There
exists a nonempty Zariski-open subset $V\subset\mathfrak{m}$ such
that for any $x\in V$ the $U$-orbit of the element $x$ intersects
$\mathcal{Z}$ at a unique point}.}

\textsc{Proof.} By Theorem \ref{Exist_of_representative} there
exists a nonempty Zariski-open subset $W$ such that for any
$x\in\mathfrak{m}$ there exists $g\in N$ satisfing
$\mathrm{Ad}_gx\in \mathcal{Y}$. Fix any $x\in W$, there is an
element $g\in N$ corresponding to $x$. Since $N=U_L\ltimes U$, $g\in
N$ can be represented as the product $g=g_1g_2$, where $g_1\in U_L$
and $g_2\in U$. Then $g_1^{-1}g\in U$. Let us show that we can take
$V=W$ and $\mathrm{Ad}_{g_1^{-1}g}x\in \mathcal{Z}$.

Since $\mathcal{Y}\subset \mathcal{Z}$ and one-parameter subgroups
$g_{i,j}(t)=I+tE_{i,j}$, where
$(i,j)\in\Delta^{\!+}_{\mathfrak{r}}$, generate the group $U_L$, it
is enough to show that for any $g_{i,j}(t)\in U_L$ we have
$\mathrm{Ad}_{g_{i,j}(t)}\mathcal{Z}\subset\mathcal{Z}$. Suppose
$g_{i,j}(t)\in U_L$; then $(i,j)\in\Delta^{\!+}_{\mathfrak{r}}$.
This means that there exists $k$ such that $R_{k-1}<i<j\leqslant
R_k$. If some element is changed after the action of the
one-parameter subgroup $g_{i,j}(t)$, then this element is $(i,a)$ or
$(b,j)$ for some $a>i$ and $b<j$. In the first case the $j$th row is
added to the $i$th row
$$\mathrm{Ad}_{g_{i,j}(t)}x_{i,a}=x_{i,a}+tx_{j,a}.$$
We have that the variables $x_{(i,a)}$ and $x_{(j,a)}$ are in the
same block $X_{k,l}$. In the second case the $i$th column is added
to the $j$th column
$$\mathrm{Ad}_{g_{i,j}(t)}x_{b,j}=x_{b,j}-tx_{b,i}.$$
Similarly, the variables $x_{(b,j)}$ and $x_{(b,i)}$ are in the same
block $X_{m,k}$. By the definition of $T$, in the case $(i,a)$ we
have that if the root $(j,a)\in T$, then $(i,a)\in T$. This means
that the $g_{i,j}$-action does not change the set
$\mathcal{Z}=\displaystyle\sum_{\xi\in T}c_{\xi}E_{\xi}$. Similarly,
if $(b,i)\in T$, then $(b,j)\in T$.

By Lemmas~\ref{Th:N_xi-invariant} and~\ref{N_independ} any
$z\in\mathcal{Z}$ such that $N_{\xi}|_z\neq0$ for any $\xi\in T$ is
a representative of some $U$-orbit.~$\Box$

\medskip
Let $\mathcal{S}$ be the set of denominators generated by invariants
$N_{\xi}$, $\xi\in T$. Denote by $K[\mathfrak{m}]^U_{\mathcal{S}}$
localization of the algebra $K[\mathfrak{m}]^U$ on $\mathcal{S}$.
Let $$\pi:K[\mathfrak{m}]^U\rightarrow K[\mathcal{Z}]$$ be the
restriction homomorphism, $f\in K[\mathfrak{m}]^U\mapsto
f|_{\mathcal{Z}}$, where the algebra $K[\mathcal{Z}]$ is a
polynomial algebra of variables $x_{\xi}$, $\xi\in T$. Extend $\pi$
to the mapping
$\widetilde{\pi}:K[\mathfrak{m}]^U_{\mathcal{S}}\rightarrow
K[c_{\xi_1}^{\pm1},c_{\xi_2}^{\pm1},\ldots,c_{\xi_s}^{\pm1}],$ where
$\xi_1,\xi_2,\ldots,\xi_s$ are all roots in $T$.

\thenv{Proposition\label{K(m)^U}}{\emph{Let
$\{\xi_1,\xi_2,\ldots,\xi_s\}$ be a collection of roots of the broad
base} $T$. \emph{The mapping
$\widetilde{\pi}:K[\mathfrak{m}]^U_{\mathcal{S}}\rightarrow
K[c_{\xi_1}^{\pm1},c_{\xi_2}^{\pm1},\ldots,c_{\xi_s}^{\pm1}]$ is an
isomorphism and} $K[\mathfrak{m}]^U_{\mathcal{S}}=
K[N_{\xi_1}^{\pm1},N_{\xi_2}^{\pm1},\ldots,N_{\xi_s}^{\pm1}]$.}

\textsc{Proof.} Let us show that $\widetilde{\pi}$ is a
monomorphism. Indeed, if $f\in K[\mathfrak{m}]^U_{\mathcal{S}}$
satisfies $\widetilde{\pi}(f)=0$, then $f|_{\mathcal{Z}}=0$. By
Proposition~\ref{Exist_of_representative_for_U},
$\mathrm{Ad}_U\mathcal{Z}$ is dense in $\mathfrak{m}$, therefore
$f(\mathfrak{m})=0$. So $f\equiv0$ and $\pi$ is a monomorphism.

To prove that $\widetilde{\pi}$ is an epimorphism, we will show that
for any $\xi\in T$ the element $c_{\xi}$ has a preimage in
$K[N_{\xi_1}^{\pm1},N_{\xi_2}^{\pm1},\ldots,N_{\xi_s}^{\pm1}]$. The
proof is by induction on the remoteness of $\xi$. Since for any root
$\xi\in M\setminus M'$ the polynomial $N_{\xi}=x_{\xi}$ is an
$U$-invariant, then $\widetilde{\pi}(N_{\xi})=c_{\xi}$ and the base
of induction is evident. Suppose for a root $\xi$ with remoteness
less than $k$ we have that $c_{\xi}$ has a preimage in
$K[N_{\xi_1}^{\pm1},N_{\xi_2}^{\pm1},\ldots,N_{\xi_s}^{\pm1}]$. Let
us show the statement for $k$. Consider a relation $\prec$ on $T$,
defined by $\varphi_1\prec\varphi_2$ whenever $i_1>i_2$ and
$j_1<j_2$, where $\varphi_1=(i_1,j_1)$ and $\varphi_2=(i_2,j_2)$.
Let $\xi\in T$ have a $k$th remoteness, then
$$N_{\xi}=x_{\xi}\prod_{\varphi\prec\xi}N_{\varphi}+b,$$
where the product is taken on all roots $\varphi\prec\xi$ such that
$\varphi\in S$ and $\varphi$ is maximal in the sense of the relation
$\prec$. For Example~\ref{Ex-2132-T} we have
$$\prod_{\varphi\prec(1,6)}N_{\varphi}=N_{(2,3)}N_{(3,4)}.$$
Note that all variables $c_{\gamma}$ in the polynomial
$\widetilde{\pi}(b)$ correspond to the roots $\gamma$ with less
remoteness than the remoteness of $\xi$. Therefore by the induction
assumption, for all these roots $\gamma$ we have that $c_{\gamma}$
has a preimage in the localization
$K[N_{\xi_1}^{\pm1},N_{\xi_2}^{\pm1},\ldots,N_{\xi_s}^{\pm1}]$.
Hence there is a function $\phi(y_1,\ldots,y_s)\in
K[y_1^{\pm1},y_2^{\pm1},\ldots,y_s^{\pm1}]$ such that
$\widetilde{\pi}(b)=
\widetilde{\pi}\big(\phi(N_{\xi_1},\ldots,N_{\xi_s})\big)$. Then
$$\widetilde{\pi}^{-1}(c_{\xi})=
\frac{N_{\xi}-\phi(N_{\xi_1},\ldots,N_{\xi_s})}%
{\displaystyle\prod_{\varphi\prec\xi}N_{\varphi}}\in
K[N_{\xi_1}^{\pm1},N_{\xi_2}^{\pm1},\ldots,N_{\xi_s}^{\pm1}].~\Box$$

\thenv{Theorem\label{Th-K[m]^U}}{\emph{The algebra of invariants
$K[\mathfrak{m}]^U$ is a polynomial algebra of} $N_{\xi}$, $\xi\in
T$.}

\textsc{Proof.} Let us show that for $L\in K[\mathfrak{m}]^U$ one
has
$$L\in K[N_{\xi_1},N_{\xi_2},\ldots,N_{\xi_s}],$$ where
$\{\xi_1,\xi_2,\ldots,\xi_s\}$ is a collection of roots of the broad
base $T$. By Proposi\-tion~\ref{K(m)^U}, there exists a polynomial
$f$ and integers $l_1,l_2,\ldots,l_k$ such that
\begin{equation}
L=\frac{f(N_{\xi_1},N_{\xi_2},\ldots,N_{\xi_k})}{\displaystyle
\prod_{i=1}^{k}N_{\xi_i}^{l_i}}.
\label{L=}\end{equation}

By the induction on the number of $N_{\xi}$ in the denominator it is
sufficient to prove that if $LN_{\xi}\in
K[N_{\xi_1},\ldots,N_{\xi_s}]$ for some $\xi\in T$ and for some
$L\in K[\mathfrak{m}]$, then $L\in K[N_{\xi_1},\ldots,N_{\xi_s}]$.

We fix a root $\xi$. Suppose $\xi=(i,j)$ and consider the case
$\xi\in M'$. If some root $\gamma$ in the broad base $T$ has the
form $(i-1,b)$ for some $b>j$, then denote $\mu_{\gamma}=(a,b)$ for
some $a>i$ such that $\mu_{\gamma}\not\in T$. If $\gamma=(a,j+1)$
for some $a<i-1$, then denote $\mu_{\gamma}=(a,b)$ for some $b<j$
such that $\mu_{\gamma}\not\in T$. For the other roots $\gamma\in T$
and in the case $\xi\not\in M'$ we have $\mu_{\gamma}=\gamma$.

The existence of this root $\mu_\gamma$ in the case
$\mu_\gamma\neq\gamma$ is explained as follows. Since $\xi\in M'$,
then $x_\xi$ is the block $X_{k,m}$ for some $k,m$ and $m>k+1$.
Evidently, the roots $(R_k,R_k+1)$ and $(R_{m-1},R_{m-1}+1)$ are
minimal in $M$ and belong to $S$. By definition of $T$, we have
$(R_k,u)\not\in T$ and $(v,R_{m-1}+1)\not\in T$ for $u\geqslant j$
and $v\leqslant i$. These roots can be chosen for $\mu_\gamma$.

\thenv{Example}{Let us take the root $\xi=(2,7)$. The symbol
$\bullet$ marks this root on the diagram. The roots
$\mu_\gamma=\gamma$ in $T$ are pointed out by the
symbol~$\boxtimes$. The single root $\mu_{\gamma}\neq\gamma$ is
marked by~$\odot$.}
\begin{center}\refstepcounter{diagram}
{\begin{tabular}{|p{0.1cm}|p{0.1cm}|p{0.1cm}|p{0.1cm}|p{0.1cm}|p{0.1cm}|p{0.1cm}|p{0.1cm}|c}
\multicolumn{2}{l}{{\small 1\quad 2\!\!}}&\multicolumn{2}{l}{{\small
3\quad 4\!\!}}&\multicolumn{2}{l}{{\small 5\quad 6\!\!}}&
\multicolumn{2}{l}{{\small 7\quad 8\!\!}}\\
\cline{1-8} \multicolumn{3}{|l|}{1}&$\!\boxtimes$&&&$\!\boxtimes$&&{\small 1}\\
\cline{4-8} \multicolumn{3}{|c|}{1}&$\!\boxtimes$&&&$\!\bullet$&$\!\boxtimes$&{\small 2}\\
\cline{4-8} \multicolumn{3}{|l|}{\qquad\ \ 1}&$\!\boxtimes$&&&&$\!\odot$&{\small 3}\\
\cline{1-8} \multicolumn{3}{|c|}{}&1&$\!\boxtimes$&&&&{\small 4}\\
\cline{4-8} \multicolumn{4}{|c|}{}&1&$\!\boxtimes$&$\!\boxtimes$&$\!\boxtimes$&{\small 5}\\
\cline{5-8} \multicolumn{5}{|c|}{}&\multicolumn{3}{|l|}{1}&{\small 6}\\
\multicolumn{5}{|c|}{}&\multicolumn{3}{|c|}{1}&{\small 7}\\
\multicolumn{5}{|c|}{}&\multicolumn{3}{|r|}{1}&{\small 8}\\
\cline{1-8} \multicolumn{8}{c}{Diagram \arabic{diagram}}\\
\end{tabular}}
\end{center}

Consider a set of matrices
$$A=\left\{\sum_{\gamma\in T}c_{\mu_\gamma}E_{\mu_\gamma},
\mbox{ where }c_{\mu_\gamma}
\mbox{ such that }N_{\gamma}|_A\neq0\mbox{ for }\gamma\neq\xi\mbox{
and }N_{\xi}|_A=0\right\}.$$

Consider a subset $\mathcal{X}=\{(N_{\xi_1}|_A,\ldots,
N_{\xi_s}|_A)\}$ of the vector space $K^{s}$. Evidently,
$\mathcal{X}\subset\mathrm{Ann}\,N_{\xi}$. Let us show that the
system of polynomials $$\{N_{\gamma}|_A,\  \gamma\neq\xi\}$$ is
algebraically independent. The proof is by induction on the number
of roots. Since we have $N_{\gamma}|_A=x_{\gamma}$ for any
$\gamma\in M\setminus M'$ and
$N_{\gamma}|_A=N_{\gamma}|_\mathcal{Z}$ for any $\gamma<\xi$, the
set $B=\{\gamma\in T:\ \gamma<\xi\}\cup (M\setminus M')$ is the base
of induction. Suppose that for a subset $B\subset T$ such that for
any root $\gamma\in B$ with the maximal remoteness and for any
$\eta\in T$ we have $\mu_{\eta}<\mu_{\gamma}$, then $\eta\in B$.
Suppose that the polynomials $N_{\gamma}|_A$, $\gamma\in B$, are
algebraically independent. Let $\gamma\not\in B\setminus\{\xi\}$ be
a root such that there is no $\eta\in T\setminus B$ such that
$\mu_{\eta}<\mu_{\gamma}$. Then $N_{\gamma}|_A=ax_{\mu_{\gamma}}+b$,
where polynomials $a,b$ depend on variables $x_{\mu_{\eta}}$,
$\eta<\gamma$. Therefore there is a single polynomial consisting the
variable $x_{\mu_{\eta}}$ in the list
$\{N_{\gamma}|_A,N_{\eta}|_A,\eta\in B\}$. Using the induction
hypothesis, $N_{\gamma}|_A$ and $N_{\eta}|_A$, where $\eta\in B$,
are algebraically independent.

Denote $\mathcal{I}_\mathcal{X}=\{\varphi\in
K[y_{\xi_1},\ldots,y_{\xi_s}]:\ \varphi(\mathcal{X})=0\}$ and
$\mathcal{I}=<y_{\xi}>$. Now let us prove that
$\mathcal{I}_\mathcal{X}=\mathcal{I}$. Obviously,
$\mathcal{I}_\mathcal{X}\supset\mathcal{I}$, hence
$\mathcal{X}\subset\mathrm{Ann}\,\mathcal{I}$. Since the dimension
of $\mathrm{Ann}\,\mathcal{I}$ is the degree of transcendence of the
algebra $K[y_{\xi_1},\ldots,y_{\xi_s}]/\mathcal{I}$ over the main
field $K$, we have
$$\dim\,\mathrm{Ann}\,\mathcal{I}=
\mathrm{degtr}_KK[y_{\xi_1},\ldots,y_{\xi_s}]/\mathcal{I}=s-1,$$
$$\dim\,\mathcal{X}=s-1.$$
Therefore, $\mathrm{Ann}\,\mathcal{I}=\overline{\mathcal{X}}$.
Suppose $g\in\mathcal{I}_\mathcal{X}$, then there exists $m\in
\mathbb{N}$ such that $g^m\in\mathcal{I}$ by the Hilbert's
Nullstellensatg. Since $\mathcal{I}$ is a prime ideal, we obtain
$g\in \mathcal{I}$. This means
$\mathcal{I}_\mathcal{X}=\mathcal{I}=<y_{\xi}>$. To conclude the
proof, it remains to note that there exists a polynomial
$p=p(y_{\xi_1},\ldots,y_{\xi_s})$ such that
$$LN_{\xi}=N_{\xi}p(N_{\xi_1},\ldots,N{\xi_s}).$$
Finally, we have $L\in K[N_{\xi_1},\ldots,N_
{\xi_s}]$.~$\Box$

By~\cite{Kh} one has

\thenv{Theorem (Khadzhiev)\label{Th-Khadzhiev}}{\emph{Let $H$ be a
connected reductive group and $U$ its maximal reductive subgroup}.
\emph{Then for any finitely generated algebraic} $H$-\emph{algebra
$A$ the algebra $A^U$ is finitely generated}.}

\thenv{Corollary}{\emph{The algebra of invariants
$K[\mathfrak{m}]^N$ is finitely generated}.}

\textsc{Proof.} By Theorem~\ref{Th-K[m]^U}, the algebra of
invariants $A=K[\mathfrak{m}]^U$ is finitely generated. Therefore
the algebra of invariants
$$A^{U_L}=K\big[K[\mathfrak{m}]^U\big]^{U_L}=K[\mathfrak{m}]^N$$
under the adjoint action of the reductive group $U_L$, where $U_L$
is the Levi subgroup of the parabolic group $P$, is finitely
generated too by the Khadzhiev's theorem.~$\Box$

\medskip
\textbf{Acknowledgments.}  \emph{I am very grateful to the referees
for their comments and recommendations which helped me to improve
the presentation. I also would like to thanks Anna Melnikov and
Alexander Panov for discussions.}

\textsc{Department of Mathematics, University of Haifa, Israil}\\
\textsc{Department of Mechanics and Mathematics, Samara University,
Russia}\\ \emph{E-mail address}: \verb"berlua@mail.ru"

\end{document}